\documentclass[ijoc,sglanonrev]{informs4}
\RequirePackage{amsmath}
\RequirePackage{lineno}
\RequirePackage{eqndefns-left} % For checking the display equation width and equation environment definitions %

\RequirePackage{tgtermes}
\RequirePackage{newtxtext}
\RequirePackage{newtxmath}
\RequirePackage{bm}
\RequirePackage{endnotes}
\RequirePackage{listings}
\RequirePackage{mdframed}
\RequirePackage{listings}

\OneAndAHalfSpacedXII % Current default line spacing
\RequirePackage{algorithm}
\RequirePackage{algpseudocode}
\RequirePackage{tikz}
\usepackage{xcolor}

\definecolor{ForestGreen}{RGB}{34,139,34}

\usepackage{natbib}
 \bibpunct[, ]{(}{)}{,}{a}{}{,}%
 \def\bibfont{\small}%

\lstdefinelanguage{Julia}%
{morekeywords={abstract,break,case,catch,const,continue,do,else,elseif,%
		end,export,false,for,function,immutable,import,importall,if,in,%
		macro,module,otherwise,quote,return,switch,true,try,type,typealias,%
		using,while},%
	sensitive=true,%
	alsoother={$},%
	morecomment=[l]\#,%
	morecomment=[n]{\#=}{=\#},%
	morestring=[s]{"}{"},%
	morestring=[m]{'}{'},%
}[keywords,comments,strings]%

\newcommand{\R}{\mathbb{R}}
\newcommand{\B}{\boldsymbol}

\definecolor{c_user}{HTML}{DFEDF6}
\definecolor{c_chat}{HTML}{FEF2D7}
\definecolor{c_prompt}{HTML}{EAE6EF}

\theorembodyfont{\normalfont}
\newmdtheoremenv[%
linecolor=gray,leftmargin=60,%
rightmargin=40,
backgroundcolor=c_prompt,%
innertopmargin=10pt,%
ntheorem,
skipabove=10pt, % Space above the theorem environment
skipbelow=10pt, % Space below the theorem environment
]{myprop}{Prompt}[section]

\theoremstyle{nonumberplain}
\newmdtheoremenv[%
style=nonumberedtheorem,
linecolor=gray,leftmargin=60,%
rightmargin=40,
backgroundcolor=c_user,%blue!5,%
innertopmargin=10pt,%
skipabove=10pt, % Space above the theorem environment
skipbelow=10pt, % Space below the theorem environment
ntheorem]{quser}{User:}[section]

\newmdtheoremenv[%
linecolor=gray,leftmargin=60,%
rightmargin=40,
backgroundcolor=c_chat,%yellow!10
innertopmargin=10pt,%
skipabove=10pt, % Space above the theorem environment
skipbelow=10pt, % Space below the theorem environment
ntheorem]{qchat}{ChatGPT:}[section]

\newmdtheoremenv[%
linecolor=gray,leftmargin=60,%
rightmargin=40,
backgroundcolor=c_chat,%
innertopmargin=10pt,%
skipabove=10pt, % Space above the theorem environment
skipbelow=10pt, % Space below the theorem environment
ntheorem]{qcode}{}[section]

\newcommand{\customvspace}[1]{\vspace{#1}}
\newlength{\mylength}
\EquationsNumberedThrough    % Default: (1), (2), ...
\TheoremsNumberedThrough     % Preferred (Theorem 1, Lemma 1, Theorem 2)
\ECRepeatTheorems  %  

\MANUSCRIPTNO{IJOC-0001-2024.00}

\begin{document}
%%%%%%%%%%%%%%%%

%\setcounter{page}{1}
%\DOI{10.1287/xxxx.0000.0000}%

\RUNAUTHOR{Bertsimas \& Margaritis}

% \RUNTITLE{Predictive Maintenance in Manufacturing}
\RUNTITLE{Robust and Adaptive Optimization under a Large Language Model Lens}

\TITLE{Robust and Adaptive Optimization under a Large Language Model Lens}

\ARTICLEAUTHORS{%
	\AUTHOR{Dimitris Bertsimas}
	\AFF{Sloan School of Management, Massachussetts Institute of Technology, \EMAIL{dbertsim@mit.edu}} %, \URL{}}
\AUTHOR{Georgios Margaritis}
\AFF{Operations Research Center, Massachussetts Institute of Technology, \EMAIL{geomar@mit.edu}}
} % end of the block

\ABSTRACT{%
In this paper, we explore the application of ChatGPT in the domain of Robust and Adaptive Robust Optimization.
We demonstrate that with appropriate prompting, ChatGPT can be used
to auto-formulate and solve simple Robust and Adaptive Optimization Problems.
We first develop specialized informational prompts tailored to the domains of Adaptive and Robust Optimization.
Then, we show that using these prompts, ChatGPT
is able to (i) formulate an optimization problem, (ii) adapt the problem so that it accounts for user-specified uncertainty, (iii) derive the computationally tractable robust counterpart of the problem and
(iv) produce executable code that solves the problem.
For simple Optimization Problems, we show that ChatGPT is able to perform these steps with little to no errors.
We also highlight some instances of erroneous execution of the steps.
Overall, we show that using in-context learning, ChatGPT has the potential to adapt to
higlhy specialized and niche domains, in which it would otherwise demonstrate very poor out-of-the-box performance.
%We then attribute those errors to the inherent weakness of Large Language Models
%in manipulating mathematical primitives and expressions.
}%
%Supplemental Material:
%Data Ethics & Reproducibility Note:

%\KEYWORDS{Stochastic programming, Decision support,Uncertainty, Disaster response, Optimization}

%\HISTORY{Received: Month DD, YYYY; Accepted: Month DD, YYYY; Published Online: Month DD, YYYY}

\maketitle
%%%%%%%%%%%%%%%%%%%%%%%%%%%%%%%%%%%%%%%%%%%%%%%%%%%%%%%%%%%%%%%%%%%%%%

% Text of your paper here
\section{Introduction}
Many real-world decision-making problems
can be modeled and solved as mathematical
optimization problems.
In fact, the field of Operations Research,
has found a lot of success in a variety of important real-world applications,
such as optimizing Air Transportation (\cite{barnhart_ApplicationsOperationsResearch_2003}),
Ride-Sharing (\cite{agatz_OptimizationDynamicRidesharing_2012}), Supply Chain Management
(\cite{jr_SupplyChainEngineering_2013}) and Financial Engineering (\cite{fabozzi_RobustPortfoliosContributions_2010}),
only to name a few.
However, despite the immense versatility and its variety of uses,
the field of mathematical optimization hasn't been as popular with the worldwide community
as fields like Machine Learning.
Part of the reason for this is that modeling and solving a real-world problem using Optimization
often requires the user to be a modeling expert, who is able
to recognize specific structures in the problem and exploit them to produce accurate and 
tractable optimization formulations.
On the other hand, the widespread availability of Machine Learning packages has made Machine Learning (ML)
widely accessible to everybody, since even non-experts can train and run ML models by writing just
a few lines of code.\par 
Hence, in this paper, we attempt to explore whether Optimization can be made more
accessible to non-experts using Large Language Models, such as ChatGPT.
However, we will focus our attention on a sub-field of Optimization, known as Robust Optimization (RO),
which often requires much more modeling expertise to be applied to a real-world problem.
%In many applications, after modeling the problem and
%formulating it using an optimization formulation,
%we can then use a variety of solvers to find a high-quality solution.
%In fact, modern solvers are equipped with efficient algorithms
%that allow them to solve a variety of problems, such as Linear
%Problems (Simplex method, \cite{nelder_SimplexMethodFunction_1965}),
%Convex problems (Interior point methods, \cite{wright_PrimalDualInteriorPointMethods_1997}),
%Semidefinite problems(\cite{vandenberghe_SemidefiniteProgramming_1996}) 
%and even Mixed Integer problems (Branch-and-Bound, \cite{land_AutomaticMethodSolving_1960}).
%\par 
The reason we are focusing on Robust Optimization is that 
although traditional optimization can provide high-quality solutions in a variety of
cases, it has one major drawback; It doesn't always
account for the uncertainty that is innate to most real-world applications.
In fact, systems optimized in the classical sense can be very sensitive to small changes.
This inherent sensitivity gives rise to Robust Optimization (RO),
the sub-filed of optimization that focuses on finding robust solutions, rather than optimal ones (\cite{beyer_RobustOptimizationComprehensive_2007}).
RO helps address (a) the uncertainty in the data of an optimization problem,
(b) the implementation errors of a feasible solution to the problem and (c) multi-stage optimization
under uncertainty (\cite{bertsimas_RobustAdaptiveOptimization_}).\par 
In traditional RO, we solve the problem only using current data and we 
make here-and-now decisions for all the optimization variables of the problem.
However, another important sub-field of RO is Adaptive Robust Optimization (ARO),
which allows us to postpone part of the decision, until we have more information
on the uncertain parameters.
In particular, ARO discerns between here-and-now decision variables and wait-and-see
decision variables, where the former are determined here and now,
while the latter are determined after the values of the uncertain parameters are revealed (\cite{bertsimas_RobustAdaptiveOptimization_}).
ARO is particularly effective in multi-period (dynamic) decision problems, and generally gives
less conservative and more realistic solution than traditional RO, although at a higher
computational cost.\par 
Both RO and ARO have been applied to a variety of real-world applications,
including Portfolio Optimization (\cite{fabozzi_RobustPortfoliosContributions_2010}),
Supply Chain Management, Circuit Design and Power Control (\cite{bertsimas_TheoryApplicationsRobust_2011}).
However, in order to introduce robustness to a classical optimization problem,
there are a variety of challenges.
We first need to identify the uncertain parameters of the
problem and model them using appropriately chosen uncertainty sets.
Then, we need to embed those uncertainties into a unified RO/ARO formulation.
Finally, by applying a range of theoretical results, 
we need to convert the original RO/ARO formulation into a precise and tractable optimization
problem, which can be solved by traditional optimizers efficiently.
This process can, in practice, be complex and time-consuming,
and needs to be done by OR/modeling experts that have an very good background in RO.\par 
Clearly, this makes the application of RO much less accessible to modelers without an adequate RO background. 
Hence, in this paper, we attempt to Large Language Models (LLMs)
to help make the RO modeling process more accessible to non-experts.
In particular, we utilize ChatGPT and in-context learning %\cite{reynolds_PromptProgrammingLarge_2021}
in order to automate (i) the translation of a classical optimization formulation
to a unified RO/ARO formulation, (ii)
the translation of the original RO/ARO formulation to a tractable optimization problem and (iii) the generation of executable code that solves the Optimization Problem.

\subsection{Related Work}
The idea of using Natural Language Processing (NLP) techniques
to help in formulating optimization problems has been studied
to some extent in the literature.
In (\cite{ramamonjison_NL4OptCompetitionFormulating_2023}), the authors
propose the competition NL4OPT, which aims to gather interest
for the development of machine learning models that turn natural language
to optimization formulations.
The authors outline some of the winning solutions to the competition,
most of which utilize either some older NLP models, such as Bidirection LSTMs and
Conditional Random Fields (\cite{huang_BidirectionalLSTMCRFModels_2015}), or some of the early transformer-based models which are variations BERT (\cite{devlin_BERTPretrainingDeep_2019}).
The authors also tested ChatGPT in the competition dataset, showing very good performance, and in some cases,
even surpassing the winning solutions in certain subtasks of the competition,
despite the fact that those solutions, unlike ChatGPT, were highly taylored to the competition.\par 
Apart from the NL4OPT competition, there have also been other independent approaches which
aim at helping in the formulation of optimization problems starting from natural language.
One such approach is (\cite{ramamonjison_AugmentingOperationsResearch_2022}), in which
the authors propose a multi-step approach to do so, which includes the use of BERT-based
models and the generative LLM BART (\cite{lewis_BARTDenoisingSequencetoSequence_2019}).
Finally, in the position paper (\cite{tsouros_HolyGrailNatural_2023}), the authors propose some
ideas and examples of how ChatGPT can be used to model simple optimization problems and auto-generate
code to solve them.\par 
Although the ChatGPT-based approaches described in the literature bear some resemblance to our approach, we take a different 
turn and focus solely on RO and ARO problems, which come with their own sets of difficulties.
In particular, since RO and particularly ARO are niche research areas, ChatGPT performs very poorly
out-of-the-box in such problems, due to the limited training data that exists for those topics.
In addition to that, formulating tractable RO and ARO problems requires the modeler
to apply and combine a variety of mathematical theorems, in order to transform the problem into a tractable
formulation that can be solved using standard optimization solvers (e.g. Gurobi, CPLEX).
%In addition to that, formulating RO and ARO problems require the modeler
%to apply and combine a variety of mathematical theorems, in order to transform the problem into a tractable
%formulation that can be solved using standard optimization solvers (e.g. Gurobi, CPLEX).
Again, due to the niche nature of RO and ARO, ChatGPT is not equipped with the mathematical tools
required to perform those transformations.\par 

%Hence, in our approach, we attempt to use in-context learning
%in order to address those difficulties which are particular to RO and ARO optimization problems.
%In particular, we design  prompting schemes which allow ChatGPT to take an Optimization Problem,
%formulate it in a way which accounts for uncertainty, transform it into a tractable, optimization problem
%and then write some example code to solve it.
%To our knowledge, the process of using ChatGPT or other LLMs to formulate and solve Robust and Adaptive Robust Optimization
%problems hasn't been studied in the past, so we argue that this is the first approach in that direction.\par 

In order to overcome this difficulties, we design specific prompting schemes which can help ChatGPT auto-formulate and solve RO and ARO problems through in-context learning.
%In our approach, we attempt to use in-context learning
%in order to address those difficulties which are particular to RO and ARO optimization problems.
%In particular, we design  prompting schemes which allow ChatGPT to take an Optimization Problem,
%formulate it in a way which accounts for uncertainty, transform it into a tractable, optimization problem
%and then write some example code to solve it.
To our knowledge, the process of using ChatGPT or other LLMs to manipulate and solve RO and ARO problems hasn't been studied in the past, so we argue that our approach is the first approach in that direction.\par 

\subsection{Contribution}
%The importance of this work is that it may allow non-expert users to take a standard optimization problem,
%add robustness to it, transform the problem into a tractable optimization problem and then solve
%it with minimal knowledge of Robust and Adaptive Robust Optimization.
%Another contribution of this work is that it shows how effective in-context learning can be, 
%even in areas where Generative LLMs have very little knowledge about.
%In this spirit, we also see the power of ChatGPT in applying and combining non-trivial mathematical theorems 
%which were never in the training set, demonstrating significant promise in that direction as well.

The contribution of our work is two-fold.
Firstly, we show that with appropriate prompting, generative models like ChatGPT can learn how to reformulate RO and ARO models.
Specifically, for basic uncertainty sets and simple Optimization Problems, we demonstrate that ChatGPT can be trained to (i) account for uncertainty in an optimization problem, (ii) derive the computationally-tractable robust counterpart of the problem and (iii) produce executable code that solves the problem.
We also showcase that when dealing with more complicated problems and uncertainty sets, ChatGPT may make errors in some of those steps.\par 
Secondly, aside from the specific area of RO, we show that Large Language Models can adapt to very niche and specific areas, to which they weren't explicitly trained on.
We show that this adaptation can happen with very little in-context training and informational prompting, a procedure which resembles the human learning process.

%\clearpage
\section{Approach}
In this chapter, we will first lay out the foundations of the RO and ARO
problems that we are trying to address.
We will describe the necessary mathematical primitives that can be used
to transform an optimization formulation to a formulation that accounts for uncertainty
and finally to a tractable optimization formulation that is amenable to standard solution techniques.
Then, we will describe for both RO and ARO how we can use ChatGPT 3.5 to automate this process
using both in-context-learning and chain-of-thought prompting techniques.
\subsection{Robust Optimization}
In this section, we will first describe the necessary theoretical foundation
required to formulate RO problems and convert them into tractable optimization
formulations.
Then, we will attempt to teach ChatGPT how to use those tools in order to add uncertainty
to Linear Optimization problem, transform the problem into a tractable formulation
and then produce a code snippet that solves the problem.
\subsubsection{Background}\label{section:bg}
For the purpose of this work, we will concern ourselves with Linear Optimization problems (LO)
of the following form:
\begin{equation} \label{eq:lo1}
	\begin{aligned}
		\min_{\bm{x}\in \R^n} \text{  } &\bm{c}^T\bm{x}\\
		\text{s.t. } &\bm{Ax}\leq \bm{b} \\
		%& \beta_i\geq 0, \quad \forall i \in [m]
	\end{aligned}
\end{equation}
where $\bm{c}\in\R^n,\:\bm{b}\in\R^m$ and $\bm{A}\in\R^m\times\R^n$.
In this setting, we assume that the data $(\bm{c},\bm{A},\bm{b})$ are uncertain,
but are known to reside in an uncertainty set $\mathcal{U}$.
Hence as described in (\cite{bertsimas_RobustAdaptiveOptimization_}), 
Problem (\ref{eq:lo1}) is not a single deterministic problem,
but its rather a family of problems, one for each realization of the uncertain
parameters $(\bm{c},\bm{A},\bm{b})$ :
\begin{equation}\label{eq:ro1}
	\min_{\bm{x}}\{\bm{c}^T\bm{x}\:|\:\bm{Ax}\leq\bm{b}\}_{(\bm{c},\bm{A},\bm{b})\in\mathcal{U}}
\end{equation}
Following (\cite{bertsimas_RobustAdaptiveOptimization_}), 
we make the assumption that the user has specified a particular uncertainty set $\mathcal{U}$ on the data $(\bm{c},\bm{A},\bm{b})$
and that the optimal solution $\bm{x}_*$ has to satisfy the 
inequality constraints $\bm{Ax}\leq \bm{b}$ for all the values $\bm{A},\bm{b}$ such that
$(\bm{c},\bm{A},\bm{b})\in\mathcal{U}$.
Under those assumptions, the RO problem (\ref{eq:ro1}) can be written
as the following deterministic optimization problem:
\begin{equation}\label{eq:ro2}
	\begin{aligned}
		\min_{\bm{x},t} \text{  } &t\\
		\text{s.t. } &\bm{c}^T\bm{x}-t\leq 0,\quad (\bm{c},\bm{A},\bm{b})\in\mathcal{U} \\
		&\bm{a}_i^T\bm{x}-b_i \leq 0,\quad (\bm{c},\bm{A},\bm{b})\in\mathcal{U} \\
	\end{aligned}
\end{equation}
The resulting problem is an LO problem but with perhaps an infinite number of constraints (one for each
realization of the uncertainty parameters).
Also notice that although in (\ref{eq:ro1}) we had uncertainty both in the cost function and the constraints, 
we now have an LO problem where the uncertainty is only in the constraints, as shown in (\ref{eq:ro2}).
This means that we only need to consider uncertainty in the constraints and not in the objective.\par 
Additionally, it turns out that instead of considering uncertainty on the entire matrix $\bm{A}$ and vector $\bm{b}$
at once, we can instead consider row-wise uncertainty on the rows of $\bm{A}$ individually, with independent constraints of the following form:
\begin{equation}
	\bm{a}^T_i\bm{x}-b_i\leq 0,\quad \forall (\bm{a}_i,b_i)\in\mathcal{U}_i
\end{equation}
where $\mathcal{U}_i$ is a projection of $\mathcal{U}$ into the uncertain parameters of the $i$-th constraint.
Finally, we can assume that $b_i$ is certain and consider uncertainty only in the $\bm{a}_i$,
since if $b_i$ is uncertain, the problem can be transformed into a new uncertain problem where only the variable
coefficients are uncertain.\par 
Hence, although the original robust problem (\ref{eq:ro1}) had uncertainty in the objective,
in $\bm{A}$ and in $\bm{b}$, in fact, we only need to consider problems with row-wise uncertainty on the coefficients
of the variable $\bm{x}$.
This means that from now on, we will only consider problems of the following form:
\begin{equation}\label{eq:ro3}
	\begin{aligned}
		\min_{\bm{x}\in\R^n} \text{  } &\bm{c}^T\bm{x}\\
		\text{s.t. } &\bm{a}^T_i\bm{x}\leq b_i,\quad \bm{a}_i\in\mathcal{U}_i \\
	\end{aligned}
\end{equation}
since all of the uncertain problems of the form shown in (\ref{eq:ro1}) can be reduced to Problem (\ref{eq:ro3}).
Finally, following (\cite{bertsimas_RobustAdaptiveOptimization_}), we rewrite Problem (\ref{eq:ro3}) as:\
\begin{equation}\label{eq:rofinal}
	\begin{aligned}
		\min_{\bm{x}\in\R^n} \text{  } &\bm{c}^T\bm{x}\\
		\text{s.t. } &(\bar{\bm{a}}_i+\B{Pz}_i)^T\bm{x}\leq b_i,\quad \bm{z}\in\mathcal{Z}_i,\: i\in[m] \\
	\end{aligned}
\end{equation}
Let's now introduce some results from (\cite{bertsimas_RobustAdaptiveOptimization_}) that will help
us reformulate Problem (\ref{eq:rofinal}) into a tractable optimization problem:\par 
The indicator function of a set $\mathcal{S}$ is denoted as:
\begin{equation}
	\delta(\bm{x}|S)=\begin{cases}
		0, \quad & \text{ if }\bm{x}\in\mathcal{S}\\
		\infty, \quad & \text{otherwise}\\
	\end{cases}
\end{equation}
Then, the conjugate function of $\delta(\bm{x}|S)$ is defined as:
\begin{equation}
	\delta^*(\bm{y}|\mathcal{S})=\sup_{\bm{x}}\{\bm{y}^T\bm{x}-\delta(\bm{x}|\mathcal{S})\}=\sup_{\bm{x}\in\mathcal{S}}\bm{y}^T\bm{x}
\end{equation}
Then, it turns out that Problem (\ref{eq:rofinal}) is equivalent to the following problem:
\begin{equation}\label{eq:ro_conj}
	\begin{aligned}
		\min_{\bm{x}\in\R^n} \text{  } &\bm{c}^T\bm{x}\\
		\text{s.t. } &\bar{\bm{a}}^T_i\bm{x}+\delta^*(\bm{P}^T\bm{x}|\mathcal{Z}_i)\leq b_i,\quad\: i\in[m] \\
	\end{aligned}
\end{equation}
which is called the Robust Counterpart (RC) of the original problem.
Depending on the choice of the uncertainty sets $\mathcal{Z}_i$, Problem (\ref{eq:ro_conj}) can become
very tractable.
Let's now focus on a single constraint
\begin{equation}
	(\bar{\bm{a}}+\B{Pz})^T\bm{x}\leq b,\quad \bm{z}\in\mathcal{Z}
\end{equation}
and let's see what the RC of the constraint is for some specific uncertainty sets.\par 
The following results are based on (\cite{bertsimas_RobustAdaptiveOptimization_}):
\begin{itemize}
	\item ${p}$\textbf{-norm uncertainty set:} If $\mathcal{Z}=\{\bm{z}\in\R^n\: :\: \|\bm{z}\|\leq \rho\}$, then $\delta^*(\bm{P}^T\bm{x}|\mathcal{Z})=\rho\|\bm{Px}\|_q$ where $1/q+1/p=1$.
	Hence, the RC of constraint (\ref{eq:ro_conj}) is:
	\begin{equation}
		\bar{\bm{a}}^T\bm{x}+\rho\|\bm{x}\|_q\leq b
	\end{equation}
	\item \textbf{Polyhedral uncertainty set}: If $\mathcal{Z}=\{\bm{z}\in\R^n\: :\: \bm{d}-\bm{Dz}\geq 0\}$ then:
	\begin{equation}
		\begin{aligned}
			\delta^*(\bm{P}^T\bm{x}|\mathcal{Z})=\min_{\bm{u}} \text{  } &\bm{d}^T\bm{u}\\
			\text{s.t. } &\bm{D}^T\bm{u}=\bm{P}^T\bm{x} \\
			&\bm{u}\geq 0
		\end{aligned}
	\end{equation}
	and thus the RC of constraint (\ref{eq:ro_conj}) is the set of constraints:
	\begin{equation}
		\begin{cases}
			\bar{\bm{a}}^T\bm{X}+\bm{d}^T\bm{u}\leq b\\
			\bm{D}^T\bm{u}=\bm{P}^T\bm{x}\\
			\bm{u}\geq 0
		\end{cases}
	\end{equation}
	Notice that we removed the minimization since we have a less than or equal constraint
	and we are just looking for a feasible $\bm{u}$ that satisfies the constraints.
	\item \textbf{Intersection of uncertainty sets}: If $\mathcal{Z}=\cap_{i\in[p]}\mathcal{Z}_i$ and some technical conditions are met, then we have:
	\begin{equation}
		\delta^*(\bm{P}^T\bm{x}|\mathcal{Z})=\min_{\bm{v}^i,i\in[p]}\Bigg\{\sum_{i\in [p]}\delta^*(\bm{v}^i|\mathcal{Z}_i)\:\Bigg|\:\sum_{i=1}^p\bm{v}^i=\bm{P}^T\bm{x})\Bigg\}
	\end{equation}
	Hence, the RC of constraint (9) is the set of constraints:
	\begin{equation}
		\begin{cases}
			\bar{\bm{a}}^T\bm{x}+\sum_i\delta^*(\bm{v}^i|\mathcal{Z}_i)\leq b\\
			\sum_i\bm{w}^i=\bm{P}^T\bm{x}
		\end{cases}
	\end{equation}
\end{itemize}

\subsubsection{ChatGPT for RO problems}
In order to formulate and solve RO problems using ChatGPT,
we will use in-context learning (i.e. Prompt Engineering) so as 
to introduce the concept of Robust Optimization and the necessary
toolkit required for such problems.
As part of our prompting scheme, we will utilize 
the chain-of-thought prompting paradigm, in which we ask the 
generative model to decompose a multi-step process into smaller,
intermediate steps.
This way, the model can allocate more computation to processes that
require more reasoning steps, which has shown to be very efficient,
especially in symbolic tasks (\cite{wei_ChainofThoughtPromptingElicits_2023}).\par 
In particular, the steps we are going to follow are the following:
\begin{enumerate}
	\item \textbf{Initial Problem Generation}: We first ask ChatGPT to generate a nominal Linear Optimization problem that doesn't have
	uncertainty in its parameters.
	Apart from generic LO problems of the form "$\text{minimize } \: \bm{c}^T\bm{x} \text{ s.t. } \bm{Ax}\leq \bm{b}$",
	ChatGPT can also generate other well-known optimization problems, 
	such as a Network Flow Problem, a Facility Location Problem,
	a Diet Problem etc.
	Additionally, the user may also specify a custom, never-seen-before problem that they may want to generate.
	
	\item \textbf{Introduction of Uncertain Parameters}: We ask ChatGPT to reformulate the problem of the previous step so that it
	also accounts for uncertainty in the problem data.
	In this step, we also specify the types of uncertainty sets that we would like to account
	for, which can perhaps be different for every constraint.
	During this step, we also provide ChatGPT with a quick context-prompt about what is Robust Optimization,
	how we should transform the constraints in order to account for uncertainty
	and what are the available uncertainty sets.
	
	\item \textbf{Formulation of the Robust Counterpart}: We ask ChatGPT to derive the robust counterpart of the problem of Step 2 in order
	to transform the problem into a tractable optimization form.
	In this step, we provide context-prompts with the value of $\delta^*(x|\mathcal{Z})$ for different uncertainty sets,
	as shown in Section \ref{section:bg}.
	We also provide information about how to handle the combinations of uncertainty sets (e.g. intersection).
	
	\item \textbf{Code Generation}: In the last part, we ask ChatGPT to generate a runnable Julia code that solves the problem using example data.
	To do this, we also included a prompt that helps ChatGPT treat $L_1$ and $L_{\infty}$ norms (i.e. by linearizing them)
	and $L_2$ norms (by adding second-order cone constraints).
\end{enumerate}

Out of the $4$ steps we described above, we use a context prompts only in steps $2,3$ and $4$,
since ChatGPT already knows how to formulate basic optimization problems (i.e. Step 1).
%For this reason, we will mainly focus on the prompts used for Steps $2$ and $3$.
Examples of such context prompts are presented below:

\begin{itemize}
	\item \textbf{Prompt for Step 2}: For this step, we use small variations of the following prompt:
	\begin{myprop}\label{pr:step2}
		Let's assume that you have a linear optimization problem.
		If you want to add robustness with respect to a parameter vector $\bm{a}$ which is the coefficient vector of the variables in the constraints, then you first have to bring the inequality constraints that involve your parameter of interest to the form $\bm{a}^T\bm{x}<=\bm{b}$. 
		Then you can add robustness by transforming the constraint to $(\bm{a}+\bm{Pz})^T\bm{x}<=\bm{b}$, for all $\bm{z}\in Z$, where Z is the desired uncertainty set and P is a predefined parameter matrix.
		Some common uncertainty sets $Z$ are the $p$-norm uncertainty set $Z=\{\bm{z}\in\R^n:\|\bm{z}\|_p\leq \rho\}$ 
		and the polyhedral uncertainty set $Z=\{\bm{z}\in\R^n:\bm{d}-\bm{Dz}\geq 0\}$.
	\end{myprop}
	\customvspace{\mylength}
	\item \textbf{Prompts for Step 3}: For this step, we use variations of $2$ prompts, a simpler one and a more complicated one.
	The simpler one only includes a small family of uncertainty sets and doesn't use support functions in deriving the robust counterparts:
	\begin{myprop}\label{pr:step3_1}
		In robust optimization, if you want to take the robust counterpart and make the problem tractable, you can use the following rules: Assume that you have a constraint of the form $(a+Pz)^Tx <= b$ for all z in Z.\\\\\noindent 
		- If Z is the euclidean ball uncertainty set with radius r (i.e. $||z||_2<=r$) then the constraint can be transformed to: $a^Tx+r||P^Tx||_2<=b$.\\\\
		- If Z is a polyhedral uncertainty set (i.e. $Z={z: Dz<=d}$) then the constraint can be transformed to the following set of constraints: $ax+d^Ty<=b, D^Ty=P^Tx,y>=0$.\\\\
		- If Z is a general p-norm ball uncertainty set with radius r (i.e. $Z={z: ||z||_p<=r}$), then the constraint can be transformed to: $a^Tx+r||P^Tx||_q<=b$ where $q=p/(p-1)$. This means that $||.||_q$ is the dual norm of $||.||_p$.\\\\
		- If Z is a 1-norm uncertainty set with radius r (i.e. $Z=\{z: ||z||_1<=r\}$), then the constraint can be transformed to $a^Tx+r||x||_\infty<=b$.\\\\
		- If Z is a infinity-norm uncertainty set with radius r(i.e. $Z=\{z: ||z||_\infty<=r\}$), then the constraint can be transformed to $a^Tx+r||x||_1<=b$.
	\end{myprop}
	\customvspace{\mylength}
	%\clearpage
	The second prompt uses support functions for the derivation of the robust counterpart and includes more uncertainty sets. 
	The benefit of this prompt is that it also includes information on how to combine uncertainty sets (e.g. intersections or Minkowski sums):
	\begin{myprop}\label{pr:step3_2}
		In robust optimization, if you want to take the robust counterpart and make the problem tractable, you can use the following rules: Assume that you have a constraint of the form $(a+Pz)^Tx <= b$ for all z in Z.
		Then the robust counterpart is $a^Tx+d^*(P^Tx|Z)<=b$.
		Here, $d^*(x|Z)$ is the conjugate of the support function of the uncertainty set Z evaluated at x.
		The following properties hold for the support function:\\\\\noindent
		-If Z is the infinity norm ball (i.e. $Z={z: ||z||_{\infty}<=r}$ then $d^*(x|Z)=r||x||_1$.\\\\
		-If Z is the euclidean ball (i.e. $Z={z: ||z||_2<=r}$) then $d^*(x|Z)=r||x||_2$.\\\\
		-If Z is the p-norm ball (i.e. $Z={z: ||z||_ p<=r}$ then $d^*(x|Z)=r||x||_q$ where $q=p/(p-1)$. This means that $||.||_q$ is the dual norm of $||.||_p$.\\\\
		-If $Z=\{z:\sum_{i=1}^L(1+z_i)ln(1+z_i)+(1-z_i)ln(1-z_i)\}$ then $d^*(x|Z)=sum_{i=1}^L\{t_i-s_i+u(e^{s_i/u-1}+ e^{-t_i/u-1}\}$.\\\\
		-If Z is the polyhedral uncertainty set (i.e. $Z=\{z:Dz<=d\}$) then $d^*(x|Z)$ should be replace with $d^Tu$ and the following 2 constraints should be added: $D^Tu=P^Tx$ and $u>=0$.\\\\
		-if Z is the intersection of $Z1$ and $Z2$, then $d^*(x|Z)$ is not explicitly representable, but can be represented if you replace $d^*(x|Z)$ with $d^*(w1|Z1)+d^*(w2|Z2)$ and add the constraint x=w1+w2 where x is the argument of $d^*(x|Z)$.\\\\
		-If Z is the minkowski sum of Z1 and Z2 then $d^*(x|Z)=d^*(x|Z1)+d^*(x|Z2)$.
	\end{myprop}
	\customvspace{\mylength}
	\clearpage
	\item \textbf{Prompt for Step 4}: For the code-writing step, we used variations of the following prompt:
	\begin{myprop}\label{pr:step4}
		- If you have a constraint of the form $g(x)+||x||_1<=b$ (or a sum of absolute values), then you need to get rid of the 1-norm by introducing auxiliary variables $t_i$, replacing the constraint with $g(x)+ sum(t_i)<=b$ and adding additional constraints $t_i>=x_i$ and $t_i >= -x_i$ for each i. This process is used to linearize the l1-norm.\\\\
		-If you have a constraint of the form $g(x)+||x||_\infty<=b$, then you need to get rid of the infinity norm. You can do that by introducing auxiliary variable t and replacing the constraint with $g(x)+t<=b$ while also adding the constraints $t>=x_i$ and $t>=-x_i$ for all i. This way you can linearize the l-infinity norm.\\\\
		-If you have a constraint of the form $g(x)+||x||_2<=b$, then you need to get rid of the 2 norm. You can do that by introducing auxiliary variable t>=0, replacing the constraint with $g(x)+t<=b$ and then adding the constraints $t^2>=sum(x_i^2)$
	\end{myprop}
\end{itemize}
\customvspace{\mylength}
%\clearpage
We will now present some examples of actual interactions with ChatGPT using the prompts and the steps we mentioned.
When presenting ChatGPT responses, we may omit parts of the response to show only the most relevant information.
\begin{itemize}
	\item \textbf{Example 1}: User-defined problem with different uncertainty per constraint:
	\begin{quser}
		Take the problem maximize $50*x1+40*x2+60*x3+30*x4$ subject to $10*x1+20*x2+30*x3+40*x4>=500$ and $2*x1+3*x2+4*x3+5*x4<=300$ and add robust optimization to it. In the coefficients of the first constraint, you should add euclidean ball uncertainty, whereas in the coefficients of the second constraint, you should add infinity-norm uncertainty. Please use the following information to guide you:\\
		\leavevmode
		\textbf{\{Prompt \ref{pr:step2}\}}
	\end{quser}
	\begin{qchat}
		\leavevmode\\
		...\\
		Let's begin by converting the inequality constraints into the desired form $a^Tx <= b$ for robust optimization.\\\\
		First Constraint: $10x_1 + 20x_2 + 30x_3 + 40x_4 >= 500$.
		We can rewrite this constraint as: $-10x_1 - 20x_2 - 30x_3 - 40x_4 <= -500$
		Let's bring it to the form $a^Tx <= b: a = [-10, -20, -30, -40], b = -500.$
		
		Second Constraint: $2x_1 + 3x_2 + 4x_3 + 5x_4 <= 300$.
		This constraint is already in the desired form: $a = [2, 3, 4, 5], b = 300$
		
		Now, we'll add robustness to these constraints using the specified uncertainty sets.
		
		Euclidean Ball Uncertainty for the First Constraint:\\
		$(a + z)^Tx <= b, \forall z \in Z$, where $Z = \{z: ||z||_2 <= r\}$
		
		The transformed constraint becomes: $(-10 + z_1)*x_1 + (-20 + z_2)*x_2 + (-30 + z_3)*x_3 + (-40 + z4)*x_4 <= -500$\\\\
		
		Infinity-Norm Uncertainty for the Second Constraint:
		$(a + z)^Tx <= b, \forall z \in Z$, where $Z = \{z: ||z||_{\infty} <= r\}$
		
		The transformed constraint becomes: $(2 + z_1)*x_1 + (3 + z_2)*x_2 + (4 + z_3)*x_3 + (5 + z_4)*x_4 <= 300$
		
		Finally, we can put together the robust optimization problem:\\
		
		Maximize:\\
		50x1 + 40x2 + 60x3 + 30x4
		
		Subject to:\\
		$(-10 + z_1)*x_1 + (-20 + z_2)*x_2 + (-30 + z_3)*x_3 + (-40 + z_4)*x_4 <= -500, \forall z : ||z||_2 <= r$\\
		$(2 + z_1)*x1 + (3 + z_2)*x_2 + (4 + z_3)*x_3 + (5 + z_4)*x_4 <= 300, \forall z: ||z||_{\infty} <= r$\\\\
		...
	\end{qchat}
	\begin{quser}
		Can you now derive the robust counterpart of the problem to make the problem tractable? You can use the following information to guide you:\\
		\leavevmode
		\textbf{\{Prompt \ref{pr:step3_1}\}}
	\end{quser}
	
	\begin{qchat}
		\leavevmode\\
		...\\
		Now, the robust counterpart of the problem becomes:
		
		Objective: Maximize 50x1 + 40x2 + 60x3 + 30x4
		
		Subject to:
		
		$(2x1 + 3x2 + 4x3 + 5x4) + r*sqrt(x1^2 + x2^2 + x3^2 + x4^2) <= 300$\par
		$(10x1 + 20x2 + 30x3 + 40x4) - r*(|x1| + |x2| + |x3| + |x4|) >= 500$\\
		...
	\end{qchat}
	\begin{quser}
		Can you now write a piece of julia code that solves the problem? Please use JuMP and Gurobi. In your implementation, note the following:\\
		\leavevmode
		\textbf{\{Prompt \ref{pr:step4}\}}
	\end{quser}
	\begin{qcode}
		\noindent
		\normalfont
		\textbf{ChatGPT}:
		\bigskip
		%\begin{minted}[breaklines, autogobble]{julia}
		\begin{lstlisting}
using JuMP
using Gurobi

# Define the optimization model
model = Model(Gurobi.Optimizer)
@variable(model, x[1:4] >= 0)  # Non-negative variables x1, x2, x3, x4
@variable(model, t[1:4] >= 0)
@variable(model, t2 >= 0)

r = 0.1  # Set the value of r for the uncertainty sets

@constraint(model, constraint1, 2*x[1] + 3*x[2] + 4*x[3] + 5*x[4] + r*t2 <= 300)

for i in 1:4
	@constraint(model, abs_constraint[i], t[i] >= x[i])
	@constraint(model, abs_constraint[i], t[i] >= -x[i])
end

@constraint(model, constraint2, 10*x[1] + 20*x[2] + 30*x[3] + 40*x[4] - r*sum(x[i] for i in 1:4) >= 500)
@constraint(model, l2_norm_constraint, t2^2 >= sum(x[i]^2 for i in 1:4))
@objective(model, Max, 50*x[1] + 40*x[2] + 60*x[3] + 30*x[4])

optimize!(model)

...
		\end{lstlisting}
	\end{qcode}
	%\begin{qchat}
	
	%\end{qchat}
	%	\begin{tcolorbox}[width=\textwidth,colback={green},title={With rounded corners},colbacktitle=yellow,coltitle=blue]    
		%		Take the problem maximize $50*x1+40*x2+60*x3+30*x4$ subject to $10*x1+20*x2+30*x3+40*x4>=500$ and 2*x1+3*x2+4*x3+5*x4<=300 and add robust optimization to it. In the coefficients of the first constraint, you should add euclidean ball uncertainty, whereas in the coefficients of the second constraint, you should add 3-norm uncertainty. Please use the following information to guide you:
		%		
		%		\{Prompt \ref{pr:step2}\}
		%	\end{tcolorbox}    
	
	\customvspace{\mylength}
	We can see that in this example, ChatGPT is able to take the given problem, add robustness with the specified
	uncertainty sets, convert the robust problem into its robust counterpart and finally write a piece of
	Julia code that solves the problem.
	Notice that in this particular example, ChatGPT managed to overcome a number of challenges:
	\begin{itemize}
		\item The first constraint given was of the form $\bm{a}^T\bm{x}\geq b$, while our prompts only
		described how to introduce robustness when the constraint is of the form $\bm{a}^T\bm{x}\leq  b$.
		Hence, we can see that ChatGPT brought the constraint into the correct form first (i.e. $\bm{a}^T\bm{x}\leq b$)
		and then proceeded to add robustness.
		\item Although we had constraints with different types of uncertainty, 
		ChatGPT treated the $2$ constraints independently and added a different type of uncertainty to each one of them.
		\item When taking the robust counterpart of the second constraint, ChatGPT took into account the direction of the inequality
		and correctly added a minus "-" sign in front of the $1$-norm term.
		\item When writing the Julia code for the problem, ChatGPT correctly introduced auxiliary variables to
		model the $2$ norm as a second-order cone constraint and linearize the 1-norm.
		
	\end{itemize}
	
	\item \textbf{Example 2}: Taking the robust counterpart of a well known problem in a single step:
	\begin{quser}
		Please formulate the diet problem. Also please add a budget constraint.
	\end{quser}
	\begin{qchat}
		The diet problem is a classic optimization problem in mathematics and operations research. It involves determining the optimal combination of food items to include in a daily diet while adhering to certain nutritional constraints and a budget constraint. The goal is to minimize or maximize a specific objective, such as cost or nutritional value, subject to these constraints.
		Formally, the diet problem can be defined as follows:\\\\\noindent
		\textbf{Objective}: Minimize or maximize a specific goal, which is typically one of the following:\\
		\begin{enumerate}
			\item \textbf{Minimize Cost}: Minimize the cost of the selected food items.
			\item  \textbf{Maximize Nutritional Value}: Maximize the nutritional value (e.g., total calories, protein, vitamins) of the selected food items.
		\end{enumerate}
		\leavevmode
		\\
		\noindent
		\textbf{Variables}:\\
		Let \(x_i\) be the amount (in grams or servings) of food item \(i\) to include in the diet.\\
		
		\noindent
		\textbf{Constraints}:\\
		\begin{enumerate}
			\item \textbf{Nutritional Constraints}: Ensure that the selected diet meets specific nutritional requirements. These constraints can include limits on calories, protein, carbohydrates, fats, vitamins, minerals, etc. Formally, for each nutrient or component \(j\):
			\[\sum_{i} (x_i \cdot \text{nutrient content of } i) \geq \text{minimum required nutrient content for } j\]
			\item \textbf{Budget Constraint}: Ensure that the total cost of the selected food items does not exceed a given budget limit. Formally:
			\[\sum_{i} (x_i \cdot \text{cost of } i) \leq \text{budget limit}\]
			\item \textbf{Non-negativity Constraint}: Ensure that the amounts of food items selected are non-negative:
			\[x_i \geq 0, \text{ for all } i\]
		\end{enumerate}
		...
		%	Solving the diet problem involves finding the values of \(x_i\) that optimize the chosen objective while satisfying the nutritional and budget constraints. This can be done using various optimization techniques, such as linear programming, integer programming, or heuristics.
		%	
		%	In practice, the diet problem is often used in diet planning, menu optimization, and food production to create cost-effective and nutritionally balanced meal plans while considering budget limitations.
	\end{qchat}
	%\clearpage
	\begin{quser}
		Add robust optimization to budget constraint. Please add 1-norm uncertainty. 
		Please work step by step, especially when taking the robust counterpart. 
		You can use the following piece of information to guide you:\\
		\noindent
		\textbf{\{Prompt \ref{pr:step2}\}} + \textbf{\{Prompt \ref{pr:step3_2}\}}
	\end{quser}
	\begin{qchat}
		%To add robustness to the budget constraint in the diet problem while considering 1-norm uncertainty, we'll follow the step-by-step process you provided. First, we'll bring the budget constraint to the form \(a^Tx \leq b\), and then we'll add robustness using the rules you mentioned.
		\leavevmode\\
		...\\\\
		\textbf{Step 1: Transform the Budget Constraint}\\
		The original budget constraint is:
		\[\sum_{i} (x_i \cdot \text{cost of } i) \leq \text{budget limit}\]
		
		We can rewrite this as:
		\[\mathbf{a}^T\mathbf{x} \leq b\]
		Where:\\
		- \(\mathbf{a}\) is a vector where each element corresponds to the cost of the respective food item.\\
		- \(\mathbf{x}\) is the decision variable vector representing the quantities of each food item.\\
		- \(b\) is the budget limit.\\
		
		\textbf{Step 2: Add Robustness}\\
		Now, we'll add robustness to the budget constraint by transforming it to \((\mathbf{a} + \mathbf{Pz})^T\mathbf{x} \leq b\) for all \(\mathbf{z}\) in an uncertainty set \(\mathcal{Z}\). Here, we'll consider the 1-norm uncertainty set, where \(\mathcal{Z}\) is the 1-norm ball defined as \(\mathcal{Z} = \{\mathbf{z} : \|\mathbf{z}\|_1 \leq r\}\).
		
		The robust counterpart becomes:
		\[\mathbf{a}^T\mathbf{x} + d^*_\mathcal{Z}(\mathbf{P}^T\mathbf{x}) \leq b\]
		
		Now, let's find \(d^*_\mathcal{Z}(\mathbf{P}^T\mathbf{x})\) for the 1-norm ball:
		\[d^*_\mathcal{Z}(\mathbf{P}^T\mathbf{x}) = r\|\mathbf{P}^T\mathbf{x}\|_\infty\]
		
		Substituting this into the robust constraint:
		\[\mathbf{a}^T\mathbf{x} + r\|\mathbf{P}^T\mathbf{x}\|_\infty \leq b\]
		
		This is the robust budget constraint for the diet problem with 1-norm uncertainty.\\
		
		\textbf{Note}: To make this problem tractable, you'll need to specify the values of \(r\) and \(\mathbf{P}\) based on your specific application and the level of robustness you desire. Additionally, this robust formulation can be solved using various optimization techniques, including linear programming.
	\end{qchat}
	\begin{quser}
		Can you now write a piece of julia code that solves the problem? Please use JuMP and Gurobi. In your implementation, note the following:\\
		\leavevmode
		\textbf{\{Prompt \ref{pr:step4}\}}
	\end{quser}
	\begin{qcode}
		\noindent
		\normalfont
		\textbf{ChatGPT}:
		\bigskip
		%\begin{minted}[breaklines,autogobble]{julia}
		\begin{lstlisting}
using JuMP
using Gurobi

# Define the parameters of the problem
num_food_items = 5  # Number of food items
budget_limit = 50.0  # Budget limit
costs = [10.0, 20.0, 15.0, 25.0, 30.0]  # Cost of each food item
nutrient_content = [200, 300, 150, 250, 350]  # Nutrient content of each food item
nutrient_requirements = [1500, 2000, 1000, 1800, 2200]  # Minimum nutrient requirements

# Create a JuMP model and set the solver to Gurobi
model = Model(Gurobi.Optimizer)

# Define decision variables
@variable(model, x[1:num_food_items] >= 0)

# Introduce an auxiliary variable for linearizing the infinity norm
@variable(model, t >= 0)

# Define the robust budget constraint
@constraint(model, sum(costs[i] * x[i] for i in 1:num_food_items) + t <= budget_limit)
for i in 1:num_food_items
	@constraint(model, t >= x[i])
	@constraint(model, t >= -x[i])
end

# Add constraints to ensure nutritional requirements are met
for j in 1:length(nutrient_requirements)
	@constraint(model, sum(nutrient_content[i] * x[i] for i in 1:num_food_items) >= nutrient_requirements[j])
end

# Set the objective to minimize total cost
@objective(model, Min, sum(costs[i] * x[i] for i in 1:num_food_items))

# Solve the robust diet problem
optimize!(model)

...
		\end{lstlisting}
		%\end{minted}
	\end{qcode}
	\customvspace{\mylength}
	In this example, we can see that ChatGPT is again able to add robustness to a problem, formulate
	the robust counterpart and produce executable code that solves the problem with example data.
	Some of the challenges in this particular example were the following:
	\begin{itemize}
		\item Instead of the first and simpler prompt of Step $3$, we used the second and more advanced prompt (i.e. Prompt (\ref{pr:step3_2})),
		which uses the support function in the derivation of the robust counterpart.
		Nevertheless, we can see that ChatGPT is still able to find the correct support function and use it to take the robust counterpart.
		\item In this example, we combined Steps 2 and 3 into a single step and asked ChatGPT to add robustness and take the robust counterpart
		in the same step.
		We did that by combining the prompt for Step $2$ and the prompt for Step $3$ into a single prompt.
		However, despite using a very long prompt and asking ChatGPT to execute multiple reasoning steps at the same time,
		it was still able to execute the steps correctly.
		\item The code-writing step required from ChatGPT to generate example data and also linearize the $L_{\infty}$ norm,
		both of which steps were executed successfully.
	\end{itemize}
	%Notice that in this example, we used the second, more advanced prompt for Step $3$, which also used the support function
	%in the definition of the robust counterparts.
	%Nevertheless, ChatGPT is still able to take the robust counterpart in just $1$ step.
	
\end{itemize}

\subsection{Adaptive Robust Optimization}
In this section, we will first give the necessary foundations of Adaptive Robust Optimization (ARO).
Then, we will proceed to see how ChatGPT can be used to convert optimization problems
to ARO problems, and then solve them using the Robust Optimization tools
we presented in the previous section.
%\clearpage
\subsubsection{Background}
The main difference between RO and ARO problems, is that RO assumes
that all optimization variables correspond to here-and-now decisions,
while in ARO, the values of a subset of a decision variables are revealed,
only after the uncertain parameters have been revealed (\cite{bertsimas_RobustAdaptiveOptimization_}).
Those variables are called wait-and-see variables and are
particularly useful in multi-period and dynamic problem settings.\par 
In that setting, the general form of an ARO problem is shown below (\cite{bertsimas_RobustAdaptiveOptimization_}):
\begin{equation}\label{eq:adapt}
	\begin{aligned}
		\min_{\bm{x},\bm{y}(\cdot)} \text{  } &\bm{c}^T\bm{x}+\max_{\bm{z}\in\mathcal{Z}}\bm{d}^T\bm{y}(\bm{z})\\
		\text{s.t. } &\bm{A}(\bm{z})\bm{x}+\bm{B}(\bm{z})\bm{y}(\bm{z})\leq \boldsymbol{\beta}(\bm{z}),\quad \forall \bm{z}\in\mathcal{Z} \\
		&\bm{x}\geq \bm{0},\: \bm{y}(\bm{z})\geq \bm{0}
		%& \beta_i\geq 0, \quad \forall i \in [m]
	\end{aligned}
\end{equation}
where $\bm{x}\in \R^{n_1}$ is the first-stage here-and-now decision that 
is made before the uncertainty parameter $\bm{z}\in\R^{L}$ is realized, $\bm{y}(\bm{z})\in\R^{n_2}$ denotes the
second-stage wait-and-see decision that can be adjusted according to the realized $\bm{z}$, and $\bm{A}(\bm{z}),\bm{B}(\bm{z})\in\R^{m\times n_2}$
are coefficient matrices linear in $\bm{z}$, and $\boldsymbol{\beta}(\bm{z})$ denotes the
possibly uncertain right-hand-side which is also linear in $\bm{z}$.
An important special case is the fixed-recourse case, which appears when $\bm{B}(\bm{z})=\bm{B}$,
which means that the coefficients of the adaptive variable $\bm{y}$ don't depend on $\bm{z}$.\par 
Let's now focus on a single adaptive constraint of the following form:
\begin{equation}
	(\bar{\bm{a}}+\bm{Pz})^T\bm{x}+\bm{d}^T\bm{y}(\bm{z}) \leq b, \quad \forall \bm{z}\in\mathcal{Z}
\end{equation}
Here we have assumed a constant RHS and fixed recourse.
Note that in (\ref{eq:adapt}), we are optimizing over both $\bm{x}$ and over functions $\bm{y}(\bm{z})$, which are called decision rules.
It turns out that optimizing over all the functions $\bm{y}(\cdot)$ leads to an NP-hard problem (\cite{bertsimas_RobustAdaptiveOptimization_}).
For this reason, in order to make the problem tractable, we restrict $\bm{y}(\bm{z})$ to a particular class
of functions.
Forcing the decision rule to a specific class of functions leads to an approximation and
doesn't guarantee optimality, but it is necessary to make the problem tractable.\par 
Hence, from now on, we will focus on the important special case where $\bm{y}(\bm{z})$ is linear in $\bm{z}$:
\begin{equation}
	\bm{y}(\bm{z})=\bm{u}+\bm{V}^T\bm{z}
\end{equation}
If we apply the linear decision rule to our constraint, we get a classic RO problem where we only have to make here-and-now decisions:
\begin{equation}
	(\bar{\bm{a}}+\bm{Pz})^T\bm{x}+\bm{d}^T(\bm{u}+\bm{V}^T\bm{z}) \leq b, \quad \forall \bm{z}\in\mathcal{Z}
\end{equation}
In this case, it can be proven that the constraint can be written in the following tractable way:
\begin{equation}
	\bar{\bm{a}}^T\bm{x}+\bm{d}^T\bm{u}+\delta^*(\bm{P}^T\bm{x}+\bm{V}^T\bm{D}|\mathcal{Z})  \leq b
\end{equation}
where $\delta^*(\bm{x}|\mathcal{Z})$ corresponds to the conjugate of the support function of the set $\mathcal{Z}$,
as presented in the Robust Optimization section.
Hence, we can use exactly the same tools as we used in the RO case in order to make the problem tractable.

\subsubsection{ChatGPT for ARO problems}
In order to formulate and solve ARO problems using ChatGPT, we will use very similar steps to the ones
we used in the RO case.
The steps are described below:

\begin{enumerate}
	\item \textbf{Initial Problem Generation}: We first ask ChatGPT to generate a nominal Linear Optimization problem that doesn't have
	uncertain or adaptive variables.
	Apart from generic LO problems of the form "$\text{min } \: \bm{c}^T\bm{x} \text{ s.t. } \bm{Ax}\leq \bm{b}$",
	ChatGPT can also generate other well-known optimization problems, 
	such as a Network Flow Problem, a Facility Location Problem,
	a Diet Problem etc.
	Additionally, the user may also specify a custom, never-seen-before problem that they may want to generate.
	
	\item \textbf{Introduction of Uncertain Parameters}: 
	We ask ChatGPT to reformulate the problem of the previous step so that it also includes adaptive variables
	and uncertain parameters.
	In this step, we also specify the types of uncertainty sets and decision rules that should be used.
	During this step, we also provide ChatGPT with a quick context-prompt about what is Adaptive Robust Optimization,
	and how we can add adaptive variables and account for uncertainty in the constraint.
	
	\item \textbf{Formulation of the Robust Counterpart}: We ask ChatGPT to derive the adaptive robust counterpart of the problem of Step 2 in order
	to transform the problem into a tractable optimization form.
	This step assumes that a linear decision rule is used in the previous steps, since otherwise, it is not necessary that we can make
	the problem tractable.
	Hence, assuming an LDR, we provide ChatGPT with context-prompts with the value of $\delta^*(x|\mathcal{Z})$ for different uncertainty sets,
	as shown in Section \ref{section:bg}.
	We also provide information about how to handle the combinations of uncertainty sets (e.g. intersection).
	
	%	\item \textbf{Code Generation}: In the last part, we ask ChatGPT to generate a runnable Julia code that solves the problem using example data.
	%	To do this, we also included a prompt that helps ChatGPT treat $L_1$ and $L_{\infty}$ norms (i.e. by linearizing them)
	%	and $L_2$ norms (by adding second-order cone constraints).
\end{enumerate}
We will now present the prompts used for the different steps:
\begin{itemize}
	\item \textbf{Prompt for Step 2}: For this step, we used small variations of the following prompts:
	
	\begin{myprop}\label{pr:aro2}
		An adaptive robust optimization problem (ARO) is a problem of the following form:\\
		\begin{equation}
			\begin{aligned}
				\min_{x,y(z)} & c^Tx+\max_{z \in Z} d^Ty(z)\\
				\text{s.t} &\text{ } A(z)x+B(z)y(z)<=b(z), \forall z \in \mathcal{Z}\\
				&{x}>=0, \bm{y}({z})>=0
			\end{aligned}
		\end{equation}
		
		Here, Z in an uncertainty set, x is a first-stage here-and-now decision that is made before the uncertainty z is realized, y(z) is the second-stage wait-and-see decision that can be adjusted according to the actual data (when z is realized), and A(z), B(z) are coefficient matrices linear in z. 
		Finally, b(z) is the (possibly uncertain) right hand side.  Note that in this problem we are optimizing over all the functions y(z) and the vector x. One way to make such problems tractable is to force y(z) to belong to a particular class of functions.  This process is called "decision rule" (DR). For instance, forcing y(z) to be a linear function of z (i.e. y(z)=Uz+v) corresponds to a linear decision rule (LDR). An important special case of ARO is when B(z) is constant and doesn't depend on z, in which case we say that we have constant recourse.
	\end{myprop}
	\customvspace{\mylength}
	\item \textbf{Prompt for Step 3}: For taking the Robust Counterpart, we used the following prompt:
	\begin{myprop}\label{pr:aro3}
		Let's now assume that we have 1 ARO constraint with fixed recourse of the following form:
		$(a+Pz)^Tx+d^Ty(z)<=b$, for all z in Z, where a, P and d are parameters. Note that a is a vector, P is a matrix, z is a vector, d is a vector and b is a scalar.
		Let's also assume that y is a linear decision rule (LDR) with $y=u+Vz$. Then, to make the problem tractable, the constraint can be written equivalently as:
		$a^Tx+d^Tu+d_*(P^Tx+V^Td|Z)<=b$. This corresponds to the robust counterpart.
		Here, $d_*(x|Z)$ is the conjugate of the support function of the uncertainty set Z evaluated at x.\\\\\noindent
		The following properties hold for the support function:\\
		- If $Z$ is the infinity norm ball (i.e. $Z=\{z: ||z||_{\infty}<=r\}$ then $d_*(x|Z)=r||x||_1$.\\\\
		-If $Z$ is the euclidean ball (i.e. $Z=\{z: ||z||_2<=r\}$) then $d_*(x|Z)=r||x||_2$.\\\\
		-If $Z = \{z: \sum_i a_i*e^{(d_i)^Tz}<=r$, then $d_*(x|Z)=\sum_i{v_ilog(v_i/(ai*u)-v_i)}$ with the constraints that $Dv=P^Tx, u>=0$ and $v>=0$.\\\\
		-If Z is the intersection of $Z_1$ and $Z_2$, then $d_*(x|Z)$ is not explicitly representable, but can be represented using the following constraints:   $d_*(x|Z)=d_*(w_1|Z_1)+d_*(w_2|Z_2)$ and $x=w1+w2$. Please don't forget to add both constraints when taking the robust counterpart.
	\end{myprop}
\end{itemize}
\customvspace{\mylength}
Let's now see some example interactions with ChatGPT using the described prompts:
\begin{itemize}
	\item \textbf{Example 1}: Different ARO decision rules:
	\begin{quser}
		Could you formulate an ARO problem with a linear decision rule for y?
		Please use the following information:\\\noindent
		\textbf{\{Prompt \ref{pr:aro2}\}}
	\end{quser}
	\begin{qchat}
		\leavevmode\\...\\
		\textbf{Model Formulation:}
		
		Minimize:
		\[
		c^Tx + \max_{z \in Z} (d^T U z + d^T v)
		\]
		
		Subject to:
		\[
		A(z)x + B(z)(Uz + v) \leq b(z), \quad \forall z \in Z
		\]
		
		\[
		x \geq 0, \quad U, \quad v 
		\]
		\noindent
		\textbf{Explanation}:\\\\
		- The objective function combines the first-stage cost \(c^Tx\) with the worst-case second-stage cost \(\max_{z \in Z} (d^T U z + d^T v)\).\\\\
		- The decision variables are \(x\) (first-stage decision), \(U\) (coefficient matrix for the linear decision rule), and \(v\) (offset vector for the linear decision rule).\\\\
		- The uncertainty set \(Z\) represents the possible values of the uncertainty parameter \(z\).\\\\
		- The constraints ensure that the combined decisions satisfy the uncertain constraints \(A(z)x + B(z)(Uz + v) \leq b(z)\) for all \(z\) in the uncertainty set \(Z\).\\\\
		- The decision variables \(x\), \(U\), and \(v\) are subject to non-negativity constraints.\\\\
		- The linear decision rule \(y(z) = Uz + v\) captures the adaptability of the second-stage decisions based on the actual realization of \(z\).\\
		%		
		%		This formulation optimizes over the first-stage decision \(x\), the coefficient matrix \(U\), and the offset vector \(v\). The linear decision rule \(y(z) = Uz + v\) allows for flexibility in adapting the second-stage decisions to the realized uncertainty \(z\), while the formulation seeks to minimize the overall cost considering the worst-case scenario from the uncertainty set.
	\end{qchat}
	%\clearpage
	\begin{quser}
		Could you now try a quadratic decision rule instead of a linear one?
	\end{quser}
	%\clearpage
	\begin{qchat}
		\leavevmode\\\\
		\textbf{Model Formulation:}\\
		Minimize:
		\[
		c^Tx + \max_{z \in Z} (d^T (Uz + v + Qz \otimes z))
		\]
		
		Subject to:
		\begin{equation}
			\begin{aligned}
				A(z)x + B(z)(Uz + v + Qz \otimes z) \leq b(z), \quad \forall z \in Z\\
				x \geq 0, \quad U , \quad v , \quad Q 
			\end{aligned}
		\end{equation}
		This formulation optimizes over the first-stage decision \(x\), the coefficient matrices \(U\), \(v\), and \(Q\). The quadratic decision rule \(y(z) = Uz + v + Qz \otimes z\) allows for even more flexibility in adapting the second-stage decisions to the realized uncertainty \(z\), while the formulation aims to minimize the overall cost considering the worst-case scenario from the uncertainty set.
	\end{qchat}
	\customvspace{\mylength}
	Notice that in this example, ChatGPT formulates an ARO problem with both a quadratic and a linear decision rule.
	The challenge here is that although we have described in the prompts what LDRs are, we haven't 
	said anything about Quadratic Decision Rules.
	Nevertheless, ChatGPT is able to extrapolate the concept of a linear decision rule to a quadratic one.
	In addition to that, it correctly states that in the quadratic decision rule, both $v, U$ and $Q$ are variables.
	\item \textbf{Example 2}: Taking the Robust Counterpart:
	
	\begin{quser}
		Can you formulate me an ARO with fixed recourse, linear decision rule, constant rhs and coefficient A that is linear in z?
		Please work step by step. 
		You can use the following piece of information:\\
		\noindent
		\textbf{\{Prompt \ref{pr:aro2}\}} + \textbf{\{Prompt \ref{pr:aro3}\}}
	\end{quser}
	\begin{qchat}
		Given the information provided, we will formulate an Adaptive Robust Optimization (ARO) problem with fixed recourse, a linear decision rule, constant right-hand side, and linear coefficients. Let's break down the problem step by step:
		
		\textbf{Problem Setup}:
		We want to solve the following ARO problem:
		%\begin{equation}
		\begin{align*}
			\text{minimize} \quad & c^Tx + \max_{z \in Z} d^Ty(z) \\
			\text{subject to} \quad & A(z)x + B(z)y(z) \leq b(z), \quad \forall z \in Z \\
			& x \geq 0, \quad y(z) \geq 0
		\end{align*}
		%\end{equation}
		
		We have one ARO constraint with fixed recourse:
		\[
		(a + Pz)^Tx + d^Ty(z) \leq b, \quad \forall z \in Z
		\]
		
		We are assuming a linear decision rule: \(y(z) = u + Vz\).
		
		\textbf{Equivalent Constraint Transformation}:\\
		Using the information provided, we can rewrite the ARO constraint with the linear decision rule as follows:
		%\[
		\begin{align*}
			(a + Pz)^Tx + d^T(u + Vz) &\leq b \\
			a^Tx + P^Txz + d^Tu + d^TVz &\leq b
		\end{align*}
		%\]
		
		\textbf{Expression for \(d^T\)}:\\
		Let's express \(d^T\) in terms of its components: \(d = [d_1, d_2, \ldots, d_n]\).
		%\[
		\begin{align*}
			a^Tx + P^Txz + \sum_{i=1}^n (d_i^Tu + d_i^TVz) &\leq b \\
			a^Tx + P^Txz + \sum_{i=1}^n d_i^Tu + \sum_{i=1}^n d_i^TVz &\leq b
		\end{align*}
		%\]
		
		\textbf{Conjugate of Support Function}:
		Recall that the conjugate of the support function \(d^*_Z(x)\) of an uncertainty set \(Z\) evaluated at \(x\) is given by:
		\[
		d^*_Z(x) = \max_{z \in Z} \{ x^Tz \}
		\]
		
		Applying this definition to our ARO constraint, and using the linear decision rule \(y(z) = u + Vz\), we have:
		%\begin{equation}
		\begin{align*}
			d^*_{Z}(P^Tx + V^Td | Z) &= \max_{z \in Z} \{ (P^Tx + V^Td)^Tz \} \\
			&= \max_{z \in Z} \{ x^TPz + d^TVz \}
		\end{align*}
		%\end{equation}
		
		\textbf{Substituting in Constraint}:
		Substituting the conjugate of the support function back into our constraint, we get:
		\[
		a^Tx + d^Tu + \max_{z \in Z} \{ x^TPz + d^TVz \} \leq b
		\]
		
		\textbf{Final Form}:
		The constraint in its final form becomes:
		\[
		a^Tx + d^Tu + \max_{z \in Z} \{ x^TPz + d^TVz \} \leq b
		\]
		
		This is the formulation of the ARO constraint with fixed recourse, linear decision rule, and linear coefficients that are linear in \(z\), along with a constant right-hand side.
	\end{qchat}
	\begin{quser}
		Can you now derive the robust counterpart assuming an infinity norm uncertainty set?
	\end{quser}
	\begin{qchat}
		%		Certainly! Let's derive the robust counterpart of the ARO constraint under the assumption of an infinity norm uncertainty set \(Z = \{z : \|z\|_{\infty} \leq r\}\), where \(r\) is a parameter representing the size of the uncertainty set.
		\leavevmode\\\\
		...\\\\
		Recall that the robust counterpart involves using the conjugate of the support function of the uncertainty set. For the infinity norm uncertainty set, the conjugate of the support function is \(d^*_Z(x) = r\|x\|_1\). Now, let's substitute this into the constraint:
		
		Starting from:
		\[
		a^Tx + d^Tu + \max_{z \in Z} \{ x^TPz + d^TVz \} \leq b
		\]
		
		Using the conjugate of the support function:
		\[
		a^Tx + d^Tu + \max_{z \in Z} \{ x^TPz + d^TVz \} \leq b
		\]
		
		Substituting \(d^*_Z(x) = r\|x\|_1\):
		\[
		a^Tx + d^Tu + \max_{z \in Z} \{ x^TPz + d^TVz \} \leq b
		\]
		\[
		a^Tx + d^Tu + r\|P^Tx + V^Td\|_1 \leq b
		\]
		
		Now, let's further simplify this expression by separating the absolute value term into its components:
		
		\[
		a^Tx + d^Tu + r\left(\sum_{i=1}^{n} |(P^Tx + V^Td)_i|\right) \leq b
		\]
		
		Here, \((P^Tx + V^Td)_i\) represents the \(i\)th component of the vector \(P^Tx + V^Td\).
		
		%		This is the robust counterpart of the ARO constraint with fixed recourse, linear decision rule, linear coefficients, and an infinity norm uncertainty set \(Z = \{z : \|z\|_{\infty} \leq r\}\):
		%		
		%		\[
		%		a^Tx + d^Tu + r\left(\sum_{i=1}^{n} |(P^Tx + V^Td)_i|\right) \leq b
		%		\]
		%		
		%		In this formulation, \(r\) represents the radius of the uncertainty set, and the absolute value terms in the constraint account for the uncertainty in each component of \(P^Tx + V^Td\).
	\end{qchat}
	\customvspace{\mylength}
	We can see that in this example, ChatGPT is able to formulate a general ARO problem,
	focus on a specific version of the problem where we have a linear decision rule, fixed recourse and constant RHS,
	and finally, take the robust counterpart of the ARO problem with respect to a specified uncertainty set.

	%This particular example has the following challenges:
	%	\begin{itemize}
		%		\item When supplying the 
		%	\end{itemize}
\end{itemize}

\section{Limitations}
Through the various prompting examples we presented, we showed that
ChatGPT, in a variety of scenarios, can take optimization problems, account for uncertainty through RO and ARO,
reformulate the resulting problems using tractable optimization formulations and in some cases, even generate code
to solve the resulting robust counterparts.
However, the presented approach comes with its own sets of limitations, some of which
we will try to present in this section:
\begin{itemize}
	\item \textbf{Lack of consistency}:\\
	One of the main limitations of the approach presented is that ChatGPT doesn't always produce consistent results.
	In particular, different interactions with ChatGPT with exactly the same prompts yields different responses,
	despite the fact that those responses are generally similar.
	In addition to that, it turns out that in some of those responses, ChatGPT may produce incorrect formulations
	or executable Julia codes, despite the fact that those mistakes often don't deviate too much from the correct results.
	In the example below, we give the same prompt as input to $2$ different instances of ChatGPT,
	which leads us to $2$ slightly different answers.
	Although both of those answers are correct, different answers to the same prompts may provide different context for
	subsequent prompts, thus amplifying the differences of later stages of the pipeline even further.
	\begin{quser}
		Take the problem minimize $3x_1-x_2$ s.t. $100x_1+x_2\leq 10$ and add robust optimization to the constraint. 
		You should use the  infinity-norm uncertainty set. You can use the following piece of information to guide you:\\\noindent
		\textbf{\{Prompt \ref{pr:step2}\}} + \textbf{\{Prompt \ref{pr:step3_2}\}}
	\end{quser}
	\begin{qchat}
		\textbf{(Response \#1)}: 
		
		...
		
		In the case of the infinity-norm ball, \( Z = \{ z : ||z||_{\infty} \leq r \} \), the conjugate of the support function is \( d^*_\infty(x|Z) = r||x||_1 \).
		
		Therefore, the robust counterpart becomes:
		\[ a^Tx + d^*_\infty(P^Tx | Z) \leq b \]
		\[ 100x_1 + x_2 + r(||P_1x_1 + P_2x_2||_1) \leq 10 \]
		...
	\end{qchat}
	\begin{qchat}
		\textbf{(Response \#2)}: 
		
		...
		
		Therefore, the robust counterpart becomes:
		
		\[a^Tx + r||P^Tx||_1 \leq b\]
		
		where \(a = \begin{bmatrix} 100 \\ 1 \end{bmatrix}\), \(b = 10\), and \(r\) is a chosen scalar representing the radius of the infinity-norm ball.
		
		...
		
	\end{qchat}
	\customvspace{\mylength}
	\item \textbf{Algebraic manipulations}:
	A specific area where ChatGPT struggles more is on performing algebraic manipulations.
	The inconsistent performance of ChatGPT in performing mathematical operations has already been studied in literature (\cite{frieder_MathematicalCapabilitiesChatGPT_2023}).
	To illustrate that in our problem, we provide an example below where ChatGPT is asked to add robust optimization to the constraint $100x_1+x_2\geq 10+x_1$.
	ChatGPT first attempts to bring the constraint into the form $a^Tx\leq b$,
	but it mistakenly rearranges the terms, finding that $a=[-99,1]$ instead of the correct $a=[-99,-1]$.
	This indicates that in certain cases, ChatGPT does not accurately perform arithmetic operations between algebraic 
	expressions.
	\begin{quser}
		Take the problem minimize $2x1-3x2$ s.t. $100x1+x2>=10+x1$ and add robust optimization to the constraint. You should use the  infinity-norm uncertainty set. You can use the following piece of information to guide you:
		\\\noindent
		\textbf{\{Prompt \ref{pr:step2}\}} + \textbf{\{Prompt \ref{pr:step3_2}\}}
	\end{quser}
	\begin{qchat}
		To add robust optimization to the given linear programming problem, we first need to rewrite the constraint in the form \( a^Tx \leq b \). The original constraint is \( 100x_1 + x_2 \geq 10 + x_1 \), which can be rearranged as \( -99x_1 + x_2 \geq 10 \). Thus, \( a = [-99, 1] \) and \( b = 10 \)\\
		...
		
	\end{qchat}
	\customvspace{\mylength}
	\item \textbf{Applying mathematical theorems}:
	Another area where ChatGPT may produce inconsistent results is in identifying which mathematical theorems should be used and then actually applying those theorems.
	For instance, in the example below, ChatGPT is asked to use an uncertainty set that is the intersection of a $2$-norm uncertainty set with an infinity-norm uncertainty set.
	In this case, ChatGPT is able to correctly calculate the support functions for the $2$ uncertainty sets.
	However, it doesn't use the correct theorem to combine the uncertainty sets.
	It uses the property for the Minkowski sum of uncertainty sets, 
	instead of the property for the intersection of uncertainty sets.
	This leads to an incorrect Robust counterpart of the form $a^Tx+r_1\|P^Tx\|_1+r_2\|P^Tx\|_2\leq b$ instead of the correct one
	$a^Tx+r_1\|w_1\|_1+r_2\|w_2\|_2\leq b$, where $w_1+w_2=P^Tx$.
	
	\begin{quser}
		Take the problem minimize $2x_1-3x_2$ s.t. $3x_1+x_2\leq =10$ and add robust optimization to the constraint. You should use the intersection of a 2-norm uncertainty set with an infinity-norm ball uncertainty set. Then take the robust counterpart. Please work step by step.You can use the following piece of information to guide you:\\
		\noindent
		\textbf{\{Prompt \ref{pr:step2}\}} + \textbf{\{Prompt \ref{pr:step3_2}\}}
		%\ref{pr:step2} + \ref{pr:step3_2}
	\end{quser}
	\begin{qchat}
		To add robust optimization to the given linear programming problem with uncertainty in the constraint, you should follow these steps:\\\\\noindent
		1. Express the constraint with parameter 'a' in the form $a^T x <= b$.
		Given constraint: $3x_1 + x_2 <= 10$
		We can rewrite this as: $a^T x <= b$, where $a = [3, 1]$ and $b = 10$.\\\\\noindent
		2. Define the two uncertainty sets:\par
		- Infinity-norm ball (Z1): $Z1 = \{z : ||z||_{\infty} <= r\}$\par
		- 2-norm uncertainty set (Z2): $Z2 = \{z : ||z||_2 <= r\}$.\\\\\noindent
		3. Construct the robust counterpart for each uncertainty set separately:\par
		a. For Z1 (Infinity-norm ball):
		$d_*(x|Z1) = r * ||x||_1$ (from the properties)
		The robust counterpart for Z1 is: $a^T x + r * ||P^T x||_1 <= b$\par
		b. For Z2 (2-norm ball):
		$d_*(x|Z2) = r * ||x||_2$ (from the properties)
		The robust counterpart for Z2 is: $a^T x + r * ||P^T x||_2 <= b$\\\\\noindent
		4. Now, you need to combine the two robust counterparts for Z1 and Z2. Since the overall uncertainty set is the intersection of these two, you can add them together:\par
		$a^T x + r1 * ||P^T x||_1 + r2 * ||P^T x||_2 <= b$ where r1 and r2 are the parameters for Z1 and Z2, respectively.\\\\\noindent
		5. However, we can't directly add these two terms because they have different norms. To handle this, you need to use the Minkowski sum property. The Minkowski sum of two sets Z1 and Z2 is a new set Z such that Z = $\{z : z = z_1 + z_2, z_1 \in Z_1, z_2 \in Z_2\}$.\\
		$d_*(x|Z) = d_*(x|Z1) + d_*(x|Z2)$.
		So, $d_*(x|Z) = r1 * ||P^T x||_1 + r2 * ||P^T x||_2$\\\\\noindent
		6. Now, you can construct the robust counterpart for the intersection of Z1 and Z2:
		
		$a^T x + d_*(x|Z) <= b$
		
		$a^T x + r1 * ||P^T x||_1 + r2 * ||P^T x||_2 <= b$
		
		This is the robust counterpart for the given linear programming problem with the intersection of a 2-norm uncertainty set and an infinity-norm ball uncertainty set.
		
		%\end{enumerate}
		
	\end{qchat}
\end{itemize}

\customvspace{\mylength}
\section{Conclusion}
%In this work we demonstrated the power of ChatGPT in efficiently formulating and solving
%Robust and Adaptive Robust Optimization problems.
%We showed that with the appropriate prompting techniques, we can taylor ChatGPT
%even to advanced research topics in which it wasn't trained on.
%Then, we can use the taylored version of ChatGPT to successfully solve problems within those domains,
%even generating executable code in some scenarios.\par 
%%%%%%%%%%%%%%%%%%%%%%%%%%%%%%%%%%%%%%%%%%%%

In this work we demonstrated the power of ChatGPT in formulating and solving
Robust and Adaptive Robust Optimization problems.
In particular, through a variety of examples, 
we showcased the promise of ChatGPT in taking optimization problems, 
introducing user-specified uncertainty, converting the resulting RO problem into a tractable
optimization formulation and finally, producing executable code that solves
the resulting robust counterpart.
What makes those results surprising is that although ChatGPT wasn't originally trained on RO and ARO problems,
it is able to adapt to those advanced scientific domains through limited prompt engineering and in-context learning.\par 
On the other hand, we also showed that in some cases, despite the specialized prompting,
ChatGPT makes small errors in the execution of the various steps.
In many cases, those errors stem from the inherent weaknesses of Large Language Models in manipulating
mathematical theorems and algebraic expressions.\par 
Nevertheless, despite its weaknesses, our work aims at helping popularize Mathematical Optimization
and specifically Robust and Adaptive Optimization to non-expert users.
%At the same time, however, we also showed a few instances where ChatGPT makes errors in the execution of those steps.
%More generally, however, the goal of this work is to show promise towards democratizing Robust and Adaptive optimization
%to non expert stakeholders.
We want to showcase that Large Language Models can potentially be used to increase the accessibility 
of such expert-dominated fields to non-expert stakeholders.
In particular, LLMs that have been adapted to such specialized Optimization fields can (i) act as a helpful
learning resource for non-experts, and (ii) help users with limited knowledge of the field easily
formulate and solve their problems without requiring much domain expertise.

Such endeavors may serve to increase the visibility of the field of Mathematical Optimization, especially among people
with limited prior exposure to the subject.
After all, the ubiquitous availability of Machine Learning packages that can be used with minimal expertise,
has helped ML soar in popularity over the recent years.
Hence, the integration of LLMs with Optimization workflows and the development of methods that produce
Optimization formulations from natural language has the potential to reduce the barrier of entry
into Mathematical Optimization, thus making the field more popular among both experts and non-experts.

%After all, this accessibility to non-experts is a trait traditionally enjoyed by the field of Machine Learning,
%which has helped ML soar in popularity, thanks to the ubiquitous availability of packages that can be used
%with minimal expertise. 
%
%This accessibility to non-experts is a characteristic that has long been enjoyed by the field
%of Machine Learning, which, unlike Mathematical Optimization, has a very low barrier of entry.
%However, despite that in this paper we explored the promising potential of ChatGPT in this direction,
%there is a lot of room for improvement and future research in integrating LLMs into Optimization Workflows.
%

%\bibliographystyle{nonumber}
\bibliographystyle{informs2014}
\bibliography{ms}

\begin{thebibliography}{15}
\providecommand{\natexlab}[1]{#1}
\providecommand{\url}[1]{\texttt{#1}}
\providecommand{\urlprefix}{URL }

\bibitem[{Agatz et~al.(2012)Agatz, Erera, Savelsbergh, \protect\BIBand{}
  Wang}]{agatz_OptimizationDynamicRidesharing_2012}
Agatz N, Erera A, Savelsbergh M, Wang X (2012) Optimization for dynamic
  ride-sharing: {{A}} review. \emph{European Journal of Operational Research}
  223(2):295--303, ISSN 0377-2217,
  \urlprefix\url{http://dx.doi.org/10.1016/j.ejor.2012.05.028}.

\bibitem[{Barnhart et~al.(2003)Barnhart, Belobaba, \protect\BIBand{}
  Odoni}]{barnhart_ApplicationsOperationsResearch_2003}
Barnhart C, Belobaba P, Odoni AR (2003) Applications of {{Operations Research}}
  in the {{Air Transport Industry}}. \emph{Transportation Science}
  37(4):368--391, ISSN 0041-1655,
  \urlprefix\url{http://dx.doi.org/10.1287/trsc.37.4.368.23276}.

\bibitem[{Bertsimas et~al.(2011)Bertsimas, Brown, \protect\BIBand{}
  Caramanis}]{bertsimas_TheoryApplicationsRobust_2011}
Bertsimas D, Brown DB, Caramanis C (2011) Theory and {{Applications}} of
  {{Robust Optimization}}. \emph{SIAM Review} 53(3):464--501, ISSN 0036-1445.

\bibitem[{Bertsimas \protect\BIBand{}
  Den~Hertog(2022)}]{bertsimas_RobustAdaptiveOptimization_}
Bertsimas D, Den~Hertog D (2022) Robust and {{Adaptive Optimization}}.
  https://www.dynamic-ideas.com/books/robust-and-adaptive-optimization.

\bibitem[{Beyer \protect\BIBand{}
  Sendhoff(2007)}]{beyer_RobustOptimizationComprehensive_2007}
Beyer HG, Sendhoff B (2007) Robust optimization -- {{A}} comprehensive survey.
  \emph{Computer Methods in Applied Mechanics and Engineering}
  196(33):3190--3218, ISSN 0045-7825,
  \urlprefix\url{http://dx.doi.org/10.1016/j.cma.2007.03.003}.

\bibitem[{Devlin et~al.(2019)Devlin, Chang, Lee, \protect\BIBand{}
  Toutanova}]{devlin_BERTPretrainingDeep_2019}
Devlin J, Chang MW, Lee K, Toutanova K (2019) {{BERT}}: {{Pre-training}} of
  {{Deep Bidirectional Transformers}} for {{Language Understanding}}.
  \urlprefix\url{http://dx.doi.org/10.48550/arXiv.1810.04805}.

\bibitem[{Donald~Warsing(2013)}]{jr_SupplyChainEngineering_2013}
Donald~Warsing ARR (2013) \emph{Supply {{Chain Engineering}}: {{Models}} and
  {{Applications}}} (Boca Raton: CRC Press), ISBN 978-0-429-10251-6,
  \urlprefix\url{http://dx.doi.org/10.1201/b13184}.

\bibitem[{Fabozzi et~al.(2010)Fabozzi, Huang, \protect\BIBand{}
  Zhou}]{fabozzi_RobustPortfoliosContributions_2010}
Fabozzi FJ, Huang D, Zhou G (2010) Robust portfolios: Contributions from
  operations research and finance. \emph{Annals of Operations Research}
  176(1):191--220, ISSN 1572-9338,
  \urlprefix\url{http://dx.doi.org/10.1007/s10479-009-0515-6}.

\bibitem[{Frieder et~al.(2023)Frieder, Pinchetti, Chevalier, Griffiths,
  Salvatori, Lukasiewicz, Petersen, \protect\BIBand{}
  Berner}]{frieder_MathematicalCapabilitiesChatGPT_2023}
Frieder S, Pinchetti L, Chevalier A, Griffiths RR, Salvatori T, Lukasiewicz T,
  Petersen PC, Berner J (2023) Mathematical {{Capabilities}} of {{ChatGPT}}.

\bibitem[{Huang et~al.(2015)Huang, Xu, \protect\BIBand{}
  Yu}]{huang_BidirectionalLSTMCRFModels_2015}
Huang Z, Xu W, Yu K (2015) Bidirectional {{LSTM-CRF Models}} for {{Sequence
  Tagging}}. \urlprefix\url{http://dx.doi.org/10.48550/arXiv.1508.01991}.

\bibitem[{Lewis et~al.(2019)Lewis, Liu, Goyal, Ghazvininejad, Mohamed, Levy,
  Stoyanov, \protect\BIBand{}
  Zettlemoyer}]{lewis_BARTDenoisingSequencetoSequence_2019}
Lewis M, Liu Y, Goyal N, Ghazvininejad M, Mohamed A, Levy O, Stoyanov V,
  Zettlemoyer L (2019) {{BART}}: {{Denoising Sequence-to-Sequence
  Pre-training}} for {{Natural Language Generation}}, {{Translation}}, and
  {{Comprehension}}.
  \urlprefix\url{http://dx.doi.org/10.48550/arXiv.1910.13461}.

\bibitem[{Ramamonjison et~al.(2022)Ramamonjison, Li, Yu, He, Rengan,
  {Banitalebi-Dehkordi}, Zhou, \protect\BIBand{}
  Zhang}]{ramamonjison_AugmentingOperationsResearch_2022}
Ramamonjison R, Li H, Yu TT, He S, Rengan V, {Banitalebi-Dehkordi} A, Zhou Z,
  Zhang Y (2022) Augmenting {{Operations Research}} with {{Auto-Formulation}}
  of {{Optimization Models}} from {{Problem Descriptions}}.

\bibitem[{Ramamonjison et~al.(2023)Ramamonjison, Yu, Li, Li, Carenini, Ghaddar,
  He, Mostajabdaveh, {Banitalebi-Dehkordi}, Zhou, \protect\BIBand{}
  Zhang}]{ramamonjison_NL4OptCompetitionFormulating_2023}
Ramamonjison R, Yu TT, Li R, Li H, Carenini G, Ghaddar B, He S, Mostajabdaveh
  M, {Banitalebi-Dehkordi} A, Zhou Z, Zhang Y (2023) {{NL4Opt Competition}}:
  {{Formulating Optimization Problems Based}} on {{Their Natural Language
  Descriptions}}.

\bibitem[{Tsouros et~al.(2023)Tsouros, Verhaeghe, Kad{\i}o{\u g}lu,
  \protect\BIBand{} Guns}]{tsouros_HolyGrailNatural_2023}
Tsouros D, Verhaeghe H, Kad{\i}o{\u g}lu S, Guns T (2023) Holy {{Grail}} 2.0:
  {{From Natural Language}} to {{Constraint Models}}.

\bibitem[{Wei et~al.(2023)Wei, Wang, Schuurmans, Bosma, Ichter, Xia, Chi, Le,
  \protect\BIBand{} Zhou}]{wei_ChainofThoughtPromptingElicits_2023}
Wei J, Wang X, Schuurmans D, Bosma M, Ichter B, Xia F, Chi E, Le Q, Zhou D
  (2023) Chain-of-{{Thought Prompting Elicits Reasoning}} in {{Large Language
  Models}}. \urlprefix\url{http://dx.doi.org/10.48550/arXiv.2201.11903}.

\end{thebibliography}

% References here (outcomment the appropriate case)

% CASE 1: BiBTeX used to constantly update the references
%   (while the paper is being written).
%\bibliographystyle{informs2014} % outcomment this and next line in Case 1
%\bibliography{<your bib file(s)>} % if more than one, comma separated

%\bibliographystyle{informs2014} % outcomment this and next line in Case 1
%\bibliography{sample} % if more than one, comma separated

% CASE 2: BiBTeX used to generate mypaper.bbl (to be further fine tuned)
%\input{mypaper.bbl} % outcomment this line in Case 2

%If you don't use BiBTex, you can manually itemize references as shown below.

%\bibliographystyle{nonumber}

%\begin{thebibliography}{3}
%\providecommand{\natexlab}[1]{#1}
%\providecommand{\url}[1]{\texttt{#1}}
%\providecommand{\urlprefix}{URL }
%
%\bibitem[{Smith(2005)}]{smith2005}
%Smith J (2005) Optimal resource allocation in humanitarian logistics.
%  \emph{Journal of Operations Research} 30(2):123--135.
%  
%\bibitem[{Jones(2010)}]{jones2010}
%Jones S (2010) Stochastic programming models for humanitarian logistics.
%  \emph{INFORMS Mathematics of Operations Research} 35(4):567--580.
%
%\bibitem[{Brown(2015)}]{brown2015}
%Brown D (2015) \emph{Introduction to Stochastic Programming} (Springer).
%
%\end{thebibliography}

%%%%%%%%%%%%%%%%%
\end{document}